\documentclass{amsart}
\usepackage{amsmath,amssymb,amsthm,amsfonts,amscd,mathrsfs,mathtools}
\usepackage[utf8]{inputenc}
\usepackage{graphicx}
\usepackage[hyphens]{url}
\usepackage{hyperref}
\usepackage{tikz-cd}

\theoremstyle{plain}
    \newtheorem{thm}{Theorem}[section]
    \newtheorem{prop}[thm]  {Proposition}

\theoremstyle{definition}
    \newtheorem{defn}[thm]  {Definition}

\newcommand{\Z}{\mathbb{Z}}
\newcommand{\C}{\mathbb{C}}
\newcommand{\R}{\mathbb{R}}
\newcommand{\bH}{\mathbb{H}}
\newcommand{\D}{\Delta}
\newcommand{\bD}{\bar\Delta}
\newcommand{\sxs}{S^2\times S^2}
\newcommand{\sj}{\sqrt{j}}
\newcommand{\Id}{\mathrm{Id}}

\begin{document}

\title{On two quotients of $\sxs$}

\author{Andrea Bianchi}

\address{Max Planck Institute for Mathematics,
Vivatsgasse 7, Bonn, Germany
}

\thanks{
This work was partially supported by the
\emph{European Research Council} under the
European Union’s Horizon 2020 research and innovation programme (grant agreement No. \texttt{772960}),
by the
\emph{Danish National Research Foundation} through the \emph{Copenhagen Centre for
Geometry and Topology} (\texttt{DNRF151}), and by the Max Planck Institute for Mathematics in Bonn.
}

\email{bianchi@mpim-bonn.mpg.de}

\date{\today}

\begin{abstract}
In this note we prove that two seemingly different smooth 4-manifolds arising as quotients of $\sxs$ by free actions of $\Z/4$ are in fact diffeomorphic, answering a question of Hambleton and Hillman.
\end{abstract}
\maketitle
\section{Introduction}
In \cite{HH} Hambleton and Hillman study 4-manifolds whose fundamental group is $\Z/4$ and whose universal cover is $\sxs$; they show that all such manifolds are homotopy equivalent to each other, and that there are 4 homeomorphism types of such manifolds, of which exactly 2 have vanishing Kirby-Siebenmann invariant (and are thus stably smoothable). It is then natural to ask for concrete examples of 4-manifolds representing the different homeomorphism types, or at least the two stably smoothable ones.

In \cite[Section 11]{HH} the authors construct two smooth, closed 4-manifolds with fundamental group $\Z/4$ and with universal cover $\sxs$. We briefly recall here the construction. The starting point is the identification of the 3-dimensional lens space $L_{4,1}$ with the unit tangent bundle of $\R P^2$. For $\ell\in \{1,5\}$ we can consider the lens space $L_{8,i}$, which admits $L_{4,1}$ as a double cover. We can then form the homotopy pushout (aka double mapping cylinder)
\begin{equation}
\label{eq:pushout}
\begin{tikzcd}
 L_{4,1}\cong UT\R P^2\ar[r]\ar[d] &\R P^2\ar[d,dashed]\\
 L_{8,\ell}\ar[r,dashed]&M_\ell,
\end{tikzcd}
\end{equation}
which can be considered as a smooth manifold in a natural way. Hambleton and Hillman ask whether $M_1$ and $M_5$ are homeomorphic or even diffeomorphic. We give a positive answer to their question.
\begin{thm}
\label{thm:main}
 The manifolds $M_1$ and $M_5$ are diffeomorphic.
\end{thm}

\section{An alternative model for the two manifolds}
As explained in \cite[Section 11]{HH}, the manifold $M_1$ can also be obtained by quotienting $\sxs$ by the action of the order-4 diffeomorphism $f_1\colon \sxs\to \sxs$ sending $(u,v)\mapsto (-v,u)$.
\begin{defn}
We denote by $T\in SO(3)$ the isometry of $\R^3$ acting on the canonical basis by $e_1\mapsto e_1,e_2\mapsto e_3, e_3\mapsto -e_2$.
We denote by $f_5\colon \sxs\to \sxs$ the diffeomorphism sending $(u,v)\mapsto(-T^2v,T^2u)$.
\end{defn}
We observe that $f_5$ has order 4, just as $f_1$; in fact we have $f_1^2=f_5^2$, as both coincide with the self map of $\sxs$ sending $(u,v)\mapsto(-u,-v)$. Our next aim is to prove the following proposition
\begin{prop}
\label{prop:identification}
 There are diffeomorphisms $M_\ell\cong \sxs/f_\ell$ for $\ell\in\{1,5\}$.
\end{prop}
This is done for $\ell=1$ in \cite{HH}, but we treat the two cases simultaneously here.
\begin{defn}
We let $\D\subset \sxs$ denote the diagonal, and let $\bD\subset \sxs$ be the antidiagonal. We also let $\perp\subset\sxs$ denote the subspace of pairs $(u,v)$ of orthogonal unit vectors in $\R^3$, and finally denote by $A$ and $O$ the open subspaces of $\sxs$ of pairs $(u,v)$ of unit vectors forming an angle which is acute (i.e. strictly between $0$ and $\pi/2$), respectively obtuse (i.e. strictly between $\pi/2$ and $\pi$).
\end{defn}
We identify $A\cong(0,1)\times\perp$ as follows: given $(u,v)\in A$, i.e. a pair of distinct vectors in $S^2$ forming an acute angle, we can renormalise the width of this angle to a number in $(0,1)\cong (0,\pi/2)$, and we can obtain a pair $(u',v')$ of orthogonal unit vectors by rotating both $u$ and $v$ in the $uv$-plane, away from each other, keeping their ``barycentre'' $u+v/|u+v|$ invariant. Similarly $O\cong (0,1)\times\perp$, identifying $(0,1)\cong(\pi/2,\pi)$ in an orientation-reversing fashion.
Observe that $f_\ell$ swaps $\D$ and $\bD$, it swaps $A$ and $O$, and it preserves $\perp$. The map $A\sqcup O \to (0,1)\times\perp$ given by folding the two identifications above is $f_\ell$-equivariant.

The quotient $\sxs/f_\ell$ can be stratified by two closed strata, namely $(\D\sqcup\bD)/f_\ell$ and $\perp/f_\ell$, and an open stratum, namely $(A\sqcup O)/f_\ell$.
We observe that the map $f_\ell$ restricts to diffeomorphisms $\D\to\bD$ and $\bD\to \D$, and $f_\ell^2\colon \D\to \D$ is the antipodal map of $\D\cong S^2$, for both $i=1,5$. Hence $(\D\sqcup\bD)/f_\ell\cong\R P^2$.
Similarly $f_\ell$ restricts to diffeomorphisms $A\to O$ and $O\to A$, and $f_\ell^2\colon A\to A$ can be identified with the map $\Id_{(0,1)}\times f_\ell^2|_{\perp}$, so that the open stratum in $\sxs/f_\ell$ looks like $(0,1)\times (\perp/f_\ell^2)$.
Finally, the second closed stratum can be identified with $\perp/f_\ell$.

The open stratum is attached to the closed strata, at its two ends, according to the quotient map $\perp/f_\ell^2\to\perp/f_\ell$, which is a double covering, and $\perp/f_\ell^2\to\R P^2$, given by sending $(u,v)\mapsto\pm \frac{u+v}{|u+v|}$.

\begin{defn}
\label{defn:S3}
Our favourite model for $S^3$ is the unit sphere in $\bH$, the body of quaternions. An element of $S^3$ has thus the form $z=a+bi+cj+dk$ with $a^2+b^2+c^2+d^2=1$. We denote by $\C_i\subset\bH$ the subspace of quaternions of the form $z=a+bi$, and similarly by $\C_j=\{a+cj\}\subset\bH$. We will consider $\bH$ both as a $\C_i$ and as a $\C_j$-vector space (of dimension 2), using left multiplication by elements in $\C_i$ and $\C_j$.
We let $S^1_i=S^3\cap\C_i$ be the unit circle in $\C_i$.

We denote by $I\subset\bH$ the real span of $i,j,k$; our favourite model of $S^2$ is the unit sphere in $I$, on which $S^3$ acts on right by conjugation.
We identify the Hopf fibration with the map $\eta\colon S^3\to S^2$ sending $w\mapsto w^{-1}iw$.

We identify the unit tangent bundle $UTS^2$ with the space of pairs $(a,b)$ of orthogonal imaginary unit quaternions, and we identify
the projection $UTS^2\to S^2$ with the map $(a,b)\mapsto a$.

We realise the lens space $L_{2,1}$ as the quotient of $S^3$ by multiplication by $-1$, and we
realise the lens space $L_{4,1}$ as the quotient of $S^3$ by left multiplication by $j$.

We let $\sj$ denote the quaternion $(1+j)/\sqrt{2}$. We realise $L_{8,1}$ as quotient of $S^3$ by left multiplication by $\sj$, and realise $L_{8,5}$ as quotient of $S^3$ by the map $z\mapsto \sj jzj^{-1}$, which is the composite of left multiplication by $\sj$ and conjugation by $j$.
\end{defn}

We justify the above definition as follows. Usually one considers $S^3$ as the unit sphere in $\C^2$, and defines the Hopf fibration as the projection map $S^3\to\C P^1=S^3/S^1$ sending $(x,y)\mapsto [x:y]$; realising $\C^2$ as $\bH=\C_i^2$ with standard basis $\{1,j\}$, i.e. letting
$x=a+bi$ and $y=c+di$, we obtain the Hopf fibration as the quotient map $S^3\to S^3/S^1_i$. Since the stabiliser of $i$ in the right action of $S^3$ on $S^2$ by conjugation is precisely $S^1_i$, we may as well identify the Hopf fibration with the map $\eta$ from the definition, and we obtain an identification $S^2\cong S^3/S^1_i$ (we consider the right action of $S^3$ on $S^2$ to get here the quotient of $S^3$ by left multiplication by $S^1_i$).

For each pair $(a,b)$ of orthogonal imaginary unit quaternions there is a unique pair of opposite element $\pm w\in S^3$ such that $w^{-1}iw=a$ and $w^{-1}kw=b$; this allows us to identify $UTS^2$ with $S^3/-1$, so that the projection $UTS^2\to S^2$ becomes the natural projection $S^3/-1\to S^3/S^1_i$.

To get $\R P^2$ and $U\R P^2$, we need to quotient $S^3/S^1_i$ and $S^3/-1$ by the analogues of the antipodal map of $S^2$ and its induced map on $UTS^2$. The latter maps are given by $a\mapsto -a$ and $(a,b)\mapsto (-a,-b)$, respectively.
Consider the left multiplication by $j$: it is a map $S^3\to S^3$ descending to self maps of both spaces $S^3/-1$ and $S^3/S^1_i$; for all $w\in S^3$ we have moreover that $(jw)^{-1}i(jw)=-(w^{-1}iw)$ and $(jw)^{-1}k(jw)=-(w^{-1}kw)$, implying that left multiplication by $j$ on $S^3/-1$ and $S^3/S^1_i$ corresponds to the antipodal map on $UTS^2$ and $S^2$.

We can now use the alternative complex coordinates $(a+cj)$ and $(bi+dk)$, i.e. consider $\bH\cong\C_j^2$ with basis $\{1,i\}$. 
We can identify $L_{4,1}$ with $S^3/j$, which by the previous discussion can be identified with $UTS^2$, and we identify $L_{8,1}$ with $S^3/\sj$.

It is left to identify $L_{8,5}$. We should act on the first coordinate $a+cj$ by multiplication by $\sj$, and on the second coordinate $bi+dk$ by multiplication by $\sj^5=-\sj$; so we should send $a+bi+cj+dk\mapsto \sj(a-bi+cj-dk)$; the latter is equal to $\sj j(a+bi+cj+dk)j^{-1}$. This completes the justification of Definition \ref{defn:S3}.

Recall now that $S^3/-1$, as a Lie group, is isomorphic to $SO(3)$; more precisely, there is an isomorphism sending $\sj\mapsto T$. We can thus obtain $L_{8,\ell}$ also as quotients of $SO(3)$ by left multiplication by $T$ (for $\ell=1$) and by the map $A\mapsto T^3AT^2$ (for $\ell=5$). We can further identify $SO(3)$ with $\perp$ as follows: for each pair $(u,v)$ of orthogonal vectors in $S^2$, there is a unique $w$ such that $(w,u,v)$ is an oriented orthonormal basis of $\R^3$, and we identify $(u,v)$ with the the element of $SO(3)$ sending $e_1\mapsto w$, $e_2\mapsto u$ and $e_3\mapsto v$. In this light, left action by $T$ (or by $T^3$) on $SO(3)\cong\perp$ corresponds to replacing a pair $(u,v)$ with the pair $(v,-u)$ (or with the pair $(-v,u)$), whereas right action by $T^2$ corresponds to replacing a pair $(u,v)$ with the pair $(T^2u,T^2v)$. It is now apparent that $L_{8,\ell}\cong \perp/f_\ell$ for $\ell\in\{1,5\}$.

This concludes the identification of the terms in the homotopy pushout diagram \eqref{eq:pushout} with $\perp/f_\ell$, $(\perp_A\sqcup\perp_O)/f_\ell$ and $(\D\sqcup\bD)/f_\ell$, where $\perp_A\subset A$ is the subspace of pairs $(u,v)$ forming an angle of $\pi/4$, and $\perp_O\subset O$ is the subspace of pairs $(u,v)$ forming an angle of $3\pi/4$ (so that we can identify $A$ with the cylinder $(0,1)\times\perp_A$, and $O\cong(0,1)\times\perp_O$). Also the legs of the pushout diagram have been identified with the gluing maps of the two ends of the cylindrical open strata on the closed strata. This completes the proof of Proposition \ref{prop:identification}.

\section{A family of actions}
Let $A\in SO(3)$; then we can define a map $f_A\colon\sxs\to\sxs$ by sending $(u,v)\mapsto(-Av,A^{-1}u)$. We observe that $f_A^2$ is the map $(u,v)\mapsto(-u,-v)$, and for $A=\Id$, respectively $A=T^2$, we obtain the maps $f_1$, respectively $f_5$.

Consider now the product $\sxs\times SO(3)$ as a trivial $\sxs$-bundle over $SO(3)$, and let $\Z/4$ act on the copy of $\sxs$ lying over $A$ by $f_A$. The quotient $(\sxs\times SO(3))/\Z/4$ is a smooth fibre bundle over $SO(3)$. The fibres over $\Id$ and $T^2$ are, respectively, $M_1$ and $M_5$. Since $SO(3)$ is connected, we conclude that $M_1$ and $M_5$ are diffeomorphic, which is the statement of Theorem \ref{thm:main}.

We conclude by observing that the elements $\Id,T^2$ have order 2 in $SO(3)$, i.e. they are equal to their own inverses, but even though they are connected by a path in $SO(3)$, they are not connected by a path \emph{through order-2 elements}. This implies that the decomposition of $\sxs$ as union of $(\D\sqcup\bD)$, $(A\sqcup O)$ and $\perp$, which is respected by the actions of $f_1$ and $f_5$, is not respected by the action of the maps $f_A$ for $A$ lying in a path in $SO(3)$ joining $\Id$ and $T^2$. In conclusion, the identification $M_1\cong M_5$ should not be expected to respect the above stratification of $\sxs$.

\section{Acknowledgments}
The study of the pair of manifolds $M_1$ and $M_5$ was motivated by the quest of explicit examples of pairs of manifolds in dimension $\ge4$ that are homotopy equivalent, non-homeomorphic, and have non-homotopy equivalent configuration spaces, i.e. higher dimensional analogues of the result in \cite{LongSal} in dimension 3.
I would like to thank Tommaso Grossi, Lorenzo Guerra, Andrea Marino and Paolo Salvatore for some useful conversations related to this problem, and Jonathan Hillman for some comments on a first version of this note, which I wrote in June 2023, and for encouraging me to circulate it publicly.

\bibliography{Bibliography.bib}{}
\bibliographystyle{plain}

\end{document}